\documentclass[twoside
]{amsart}%
  \usepackage{a4wide}
  \usepackage{amsmath}
  \usepackage{amssymb}
  \usepackage{amsthm}
  \usepackage{mathrsfs}
  \usepackage{bbm}
  \usepackage{pifont}
  \usepackage{color}
  \usepackage{graphicx}

	\def\AA{{\ifmmode{\mathbbm{A}}\else{$\mathbbm{A}$}\fi}}
	\def\BB{{\ifmmode{\mathbbm{B}}\else{$\mathbbm{B}$}\fi}}
	\def\CC{{\ifmmode{\mathbbm{C}}\else{$\mathbbm{C}$}\fi}}
	\def\EE{{\ifmmode{\mathbbm{E}}\else{$\mathbbm{E}$}\fi}}
	\def\FF{{\ifmmode{\mathbbm{F}}\else{$\mathbbm{F}$}\fi}}
	\def\HH{{\ifmmode{\mathbbm{H}}\else{$\mathbbm{H}$}\fi}}
	\def\KK{{\ifmmode{\mathbbm{K}}\else{$\mathbbm{K}$}\fi}}
	\def\NN{{\ifmmode{\mathbbm{N}}\else{$\mathbbm{N}$}\fi}}
	\def\PP{{\ifmmode{\mathbbm{P}}\else{$\mathbbm{P}$}\fi}}
	\def\QQ{{\ifmmode{\mathbbm{Q}}\else{$\mathbbm{Q}$}\fi}}
	\def\RR{{\ifmmode{\mathbbm{R}}\else{$\mathbbm{R}$}\fi}}
	\def\TT{{\ifmmode{\mathbbm{T}}\else{$\mathbbm{T}$}\fi}}
	\def\UU{{\ifmmode{\mathbbm{U}}\else{$\mathbbm{U}$}\fi}}
	\def\ZZ{{\ifmmode{\mathbbm{Z}}\else{$\mathbbm{Z}$}\fi}}
	
	\def\A{{\ifmmode{\mathscr{A}}\else{$\mathscr{A}$}\fi}}
	\def\B{{\ifmmode{\mathscr{B}}\else{$\mathscr{B}$}\fi}}
	\def\C{{\ifmmode{\mathscr{C}}\else{$\mathscr{C}$}\fi}}
	\def\D{{\ifmmode{\mathscr{D}}\else{$\mathscr{D}$}\fi}}
	\def\E{{\ifmmode{\mathscr{E}}\else{$\mathscr{E}$}\fi}}
	\def\F{{\ifmmode{\mathscr{F}}\else{$\mathscr{F}$}\fi}}
	\def\G{{\ifmmode{\mathscr{G}}\else{$\mathscr{G}$}\fi}}
	\def\H{{\ifmmode{\mathscr{H}}\else{$\mathscr{H}$}\fi}}
	\def\I{{\ifmmode{\mathscr{I}}\else{$\mathscr{I}$}\fi}}
	\def\J{{\ifmmode{\mathscr{J}}\else{$\mathscr{J}$}\fi}}
	\def\K{{\ifmmode{\mathscr{K}}\else{$\mathscr{K}$}\fi}}
	\def\L{{\ifmmode{\mathscr{L}}\else{$\mathscr{L}$}\fi}}
	\def\M{{\ifmmode{\mathscr{M}}\else{$\mathscr{M}$}\fi}}
	\def\N{{\ifmmode{\mathscr{N}}\else{$\mathscr{N}$}\fi}}
	\def\O{{\ifmmode{\mathscr{O}}\else{$\mathscr{O}$}\fi}}
	\def\P{{\ifmmode{\mathscr{P}}\else{$\mathscr{P}$}\fi}}
	\def\Q{{\ifmmode{\mathscr{Q}}\else{$\mathscr{Q}$}\fi}}
	\def\R{{\ifmmode{\mathscr{R}}\else{$\mathscr{R}$}\fi}}
	\def\S{{\ifmmode{\mathscr{S}}\else{$\mathscr{S}$}\fi}}
	\def\T{{\ifmmode{\mathscr{T}}\else{$\mathscr{T}$}\fi}}
	\def\U{{\ifmmode{\mathscr{U}}\else{$\mathscr{U}$}\fi}}
	\def\V{{\ifmmode{\mathscr{V}}\else{$\mathscr{V}$}\fi}}
	\def\W{{\ifmmode{\mathscr{W}}\else{$\mathscr{W}$}\fi}}
	\def\X{{\ifmmode{\mathscr{X}}\else{$\mathscr{X}$}\fi}}
	\def\Y{{\ifmmode{\mathscr{Y}}\else{$\mathscr{Y}$}\fi}}
	\def\Z{{\ifmmode{\mathscr{Z}}\else{$\mathscr{Z}$}\fi}}

	\newtheoremstyle{slanted}
	{}
	{}
	{\slshape}
	{}
	{\bfseries}
	{.}
	{ }
	{}
	
	\theoremstyle{slanted}
	\newtheorem{theo}{Theorem}[section]
	\newtheorem{prop}[theo]{Proposition}
	
	\newtheorem{question}[theo]{Question}
	
	\newtheorem{lemma}[theo]{Lemma}
	\newtheorem{definition}[theo]{Definition}

	\def\ind#1{\mathbbmss{1}_{#1}}
	\def\egdef{:=}

\newcommand{\tend}[2]{\xrightarrow[#1\to#2]{}}
	
\def\esp#1{{\EE}\left[#1\right]}
\def\espc#1#2{{\EE}\left[#1\,\left|\,#2\right.\right]}
\def\ind#1{\mathbbmss{1}_{#1}}

\title[2-fold and 3-fold mixing]{2-fold and 3-fold mixing: why 3-dot-type counterexamples are impossible in one dimension}
\author{Thierry de la Rue}
\address{Laboratoire de Math\'ematiques Rapha\"el Salem\\
	UMR 6085 CNRS -- Universit\'e de Rouen\\
	Avenue de l'Universit\'e\\
	B.P. 12\\
	F76801 Saint-\'Etienne-du-Rouvray Cedex}
\email{thierry.de-la-rue@univ-rouen.fr}

\begin{document}
\bibliographystyle{amsplain}

\maketitle

\begin{abstract}
V.A. Rohlin asked in 1949 whether 2-fold mixing implies 3-fold mixing for a stationary process $(\xi_i)_{i\in\ZZ}$, and the question remains open today. In 1978, F. Ledrappier exhibited a counterexample to the 2-fold mixing implies 3-fold mixing problem, the so-called \emph{3-dot system}, but in the context of stationary random fields indexed by $\ZZ^2$. 

In this work, we first present an attempt to adapt Ledrappier's construction to the one-dimensional case, which finally leads to a stationary process which is 2-fold but not 3-fold mixing \emph{conditionally to the $\sigma$-algebra generated by some factor process}. Then, using arguments coming from the theory of joinings, we will give some strong obstacles proving that Ledrappier's counterexample can not be fully adapted to one-dimensional stationary processes.
\end{abstract}

\section{Introduction: Rohlin's multifold mixing problem and Ledrappier's two-dimensional counterexample}

The following work is based on two recent results concerning Rohlin's multifold mixing problem which are contained in \cite{ruedl13} and \cite{ruedl14}. It seemed to me interesting to put these results together and show them in a different light, emphasizing mainly on the underlying ideas rather than on technical details.

The object of our study is a stochastic process, that is to say a family $\xi=(\xi_i)_{i\in\ZZ}$ of random variables indexed by the set of integers, and we will always assume that these random variables take their values in a finite alphabet $\AA$. If two integers $i\le j$ are given, we will denote by $\xi_i^j$ the finite sequence $(\xi_i,\xi_{i+1},\ldots,\xi_j)$. Obvious generalization of this notation to the case where $i=-\infty$ or $j=+\infty$ will also be used.

We are more particularly interested in the case where the stochastic process is \emph{stationary}, which means that the probability of observing a given cylindrical event $E$ (\textit{i.e.} an event depending only on finitely many coordinates) at any position $i\in\ZZ$ does not depend on $i$:
\begin{equation}
\label{stationarity}
\forall \ell\ge0,\ \forall E\subset\AA^{\ell+1},\ \forall i\in\ZZ,\ 
\PP\Bigl( \xi_i^{i+\ell}\in E\Bigr) \ =\ \PP\Bigl( \xi_0^\ell\in E\Bigr).
\end{equation}
Another way to characterize the stationarity of the process is to say that its distribution is invariant by the coordinate shift: Let $T:\,\AA^\ZZ\to\AA^\ZZ$ be the transformation defined by $T(\xi)=\tilde\xi$, where for all $i\in\ZZ$, $\tilde\xi_i\egdef\xi_{i+1}$. Then the stochastic process $\xi$ is stationary if and only if the distribution of $T(\xi)$ is the same as the distribution of $\xi$.

The stochastic process $\xi$ is said to be \emph{mixing} if, considering two windows of arbitrarily large size $\ell$, what happens in one window is asymptotically independent of what happens in the second window when the distance between them tends to infinity: 
\begin{equation}
\label{2-fold mixing}
\forall \ell\ge0,\ \forall E_1,E_2\subset\AA^{\ell+1},\ 
\PP\Bigl(\xi_0^\ell\in E_1, \xi_p^{p+\ell}\in E_2 \Bigr)
- \PP\Bigl( \xi_0^\ell\in E_1 \Bigr) \PP\Bigl( \xi_p^{p+\ell}\in E_2 \Bigr)\ \tend{p}{\infty}\ 0.
\end{equation}

\subsection{Rohlin's question}

In 1949, V.A. Rohlin \cite{rohli6} proposed a strengthening of the previous definition involving more than two windows: $\xi$ is said to be \emph{3-fold mixing} if
\begin{multline}
\label{3-fold mixing}
\forall \ell\ge0,\ \forall E_1,E_2,E_3\subset\AA^{\ell+1},\\
\PP\Bigl( \xi_0^\ell\in E_1, \xi_p^{p+\ell}\in E_2 ,\xi_{p+q}^{p+q+\ell}\in E_3 \Bigr)
- \PP\Bigl( \xi_0^\ell\in E_1 \Bigr) \PP\Bigl( \xi_p^{p+\ell}\in E_2 \Bigr) \PP\Bigl( \xi_{p+q}^{p+q+\ell}\in E_3\Bigr)\ \tend{p,q}{\infty}\ 0.
\end{multline}
A straightforward generalization to $k$ windows naturally gives rise to the property of being \emph{$k$-fold mixing}. To avoid any confusion, we will henceforth call the classical mixing property defined by~\eqref{2-fold mixing}: \emph{2-fold mixing}\footnote{We must point out that in Rohlin's article, the definition of $k$-fold mixing originally involved $k+1$ windows, thus the classical mixing property was called \emph{1-fold mixing}. However it seems that the convention we adopt here is used by most authors, and we find it more coherent when translated in the language of multifold self-joinings (see section~\ref{joining-section}).}.

Rohlin asked whether any stationary process which is 2-fold mixing is also 3-fold mixing. This question is still open today, but a large number of mathematical works have been devoted to the subject. Many of these works show that 2-fold mixing implies 3-fold mixing for special classes of stationary processes (see \textit{e.g.} \cite{leono2} and \cite{totok1} for Gaussian processes, \cite{host1} for processes with singular spectrum, \cite{kalik2} and \cite{ryzhi2} for finite-rank processes).

\subsection{Ledrappier's counterexample in 2 dimensions: the 3-dot system}

In the opposite direction, Ledrappier \cite{ledra2} produced in 1978 a counterexample showing that in the case of stationary processes indexed by $\ZZ^2$ (we should rather speak of stationary \emph{random fields} in this context), 2-fold mixing does not necessarily imply 3-fold mixing. Here is a description of his example: Consider 
$$ G\ \egdef\ \left\{ \left(\xi_{i,j}\right)\in\{0,1\}^{\ZZ^2}: \forall (i,j),\ \xi_{i,j}+\xi_{i+1,j}+\xi_{i,j+1}=0 \mod 2 \right\}.
$$
Let us describe a probability law $\mu$ on $G$ by the way we pick a random element in $G$: First, use independent unbiased coin tosses  to choose the $\xi_{i,0}$ on the horizontal axis (one coin toss for each $i\in\ZZ$: these random variables are independent). Now, note that the ``3-dot rule" $\xi_{i,j}+\xi_{i+1,j}+\xi_{i,j+1}=0 \mod 2$ for each $(i,j)$ completely determines the coordinates $\xi_{i,j}$ on the upper-half plane $j\ge0$. It remains to choose the $\xi_{i,j}$ for $j<0$. For this, observe that we have yet no constraint on $\xi_{0,-1}$. We choose it with an unbiased coin toss, and then the entire line $\xi_{-1,j}$ is completely determined by the 3-dot rule. To complete the whole plane, we just have to pick each of the $\xi_{0,j}$ ($j<-1$) with a coin toss, and then inductively fill each horizontal line with the 3-dot rule. 

\begin{figure}[h]
\centering
\includegraphics{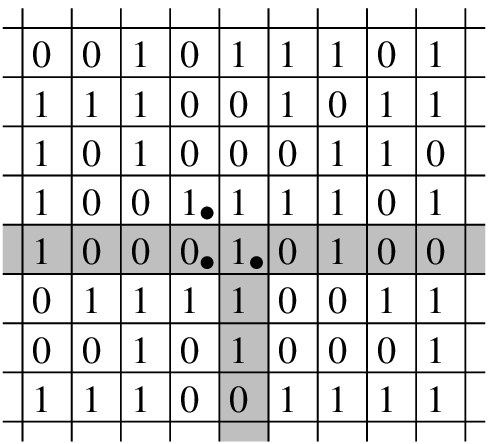}
\caption{Generation of a random configuration in $G$. First, use independent coin tosses  to choose the values of the shaded cells, then apply the 3-dot rule to complete the others: Three adjacent cells disposed as the three dotted ones must contain an even number of 1's. }
\label{fig:HaarOnG}
\end{figure}

The addition mod 2 on each coordinate turns $G$ into a compact Abelian group. We let the reader check that the probability law $\mu$ defined above on $G$ is invariant by addition of an arbitrary element of $G$, thus $\mu$ is the unique normalized Haar measure on $G$. Since any shift of coordinates in $\ZZ^2$ is an automorphism of the group $G$, such a shift preserves $\mu$. Hence $\mu$ turns $(\xi_{i,j})$ into a stationary random field. 

\begin{figure}[h]
\centering
\input{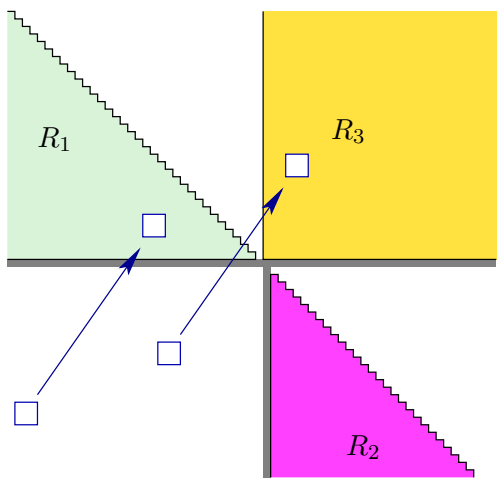}
\caption{2-fold mixing for the 3-dot system: If the distance between them is large enough, the two square windows can be shifted in such a way that one lies in one of the three colored regions, and the other one in another, independent, region.}
\label{fig:independentRegions}
\end{figure}

The definition of $k$-fold mixing for a stationary random field is formally the same as in the case of processes, except that a window is no longer an interval on the line but a square in the plane: $\{(i_0+i,j_0+j):\ 0\le i,j \le \ell\}$ for some $(i_0,j_0)\in\ZZ^2$ and some $\ell\ge0$. Let us sketch a geometric argument showing why the 2-fold mixing property holds for $(\xi_{i,j})$. Starting with the cells on the horizontal axis and the lower-half vertical axis filled with independent coin tosses , we observe that, when filling the other cells using the 3-dot rule, 
\begin{itemize}
\item the region $R_1\egdef\{(i,j):\ i<0,\ 0<j<-i\}$ only depends on the cells $(i,0)$, for $i<0$;
\item the region $R_2\egdef\{(i,j):\ j<0,\ 0<i<-j\}$ only depends on the cells $(0,j)$, for $j<0$;
\item the region $R_3\egdef\{(i,j):\ 0<i,\ 0<j\}$ only depends on the cells $(i,0)$, for $i\ge0$. (See Figure~\ref{fig:independentRegions}.)
\end{itemize}
These three regions are therefore independent. Now, if we take two windows of size $\ell$, and if the distance between them is large enough (``large enough" depending on $\ell$), it is always possible to shift the coordinates in such a way that each of the shifted windows entirely lies in one of these three regions, and not both in the same region. The two shifted windows are then independent, and since $\mu$ is preserved by coordinate shift, this means that the two windows we started with are also independent. 

\begin{figure}[h]
\centering
\input{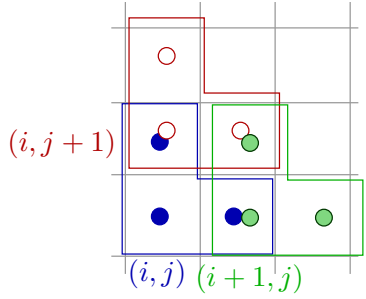}
\caption{Applying the 3-dot rule from three different corners $(i,j)$, $(i+1,j)$ and $(i,j+1)$, and adding the three equalities gives the scale-2 3-dot rule.}
\label{fig:not3foldMixing}
\end{figure}

It remains to see why Ledrappier's example is not 3-fold mixing. For this, apply the 3-dot rule from corner $(i,j)$, from corner $(i+1,j)$ and from corner $(i,j+1)$, then add the three equalities (see Figure~\ref{fig:not3foldMixing}). In the sum, the random variables $\xi_{i+1,j}$, $\xi_{i+1,j+1}$ and $\xi_{i,j+1}$ are counted twice, thus they vanish since we work modulo~2. We get the following equality, which could be called the \emph{scale-2 3-dot rule}:
\begin{equation}
\label{3-dot-bis}
\xi_{i,j}+\xi_{i+2,j}+\xi_{i,j+2}\ =\ 0\mod 2.
\end{equation}
A straightforward induction then shows that for any $n\ge0$, the \emph{scale-$2^n$ 3-dot rule} holds:
\begin{equation}
\label{3-dot-n}
\xi_{i,j}+\xi_{i+2^n,j}+\xi_{i,j+2^n}\ =\ 0\mod 2.
\end{equation}
But this shows that the three windows of size 1 $\{(i,j)\}$, $\{(i+2^n,j)\}$ and $\{(i,j+2^n)\}$ are always far from being independent, although the distance between them can be made arbitrarily large. Hence the random field $\xi$ is not 3-fold mixing.

\section{Attempt to construct a 3-dot-type one-dimensional process}

\subsection{Block construction of a 2-fold but not 3-fold mixing process}
\label{non-stationary}
In this section we describe a naive attempt to mimic the 3-dot construction on a one-dimensional process.
Our process will take its values in the same alphabet $\AA=\{0,1\}$ as for Ledrappier's example, and we start by randomly picking the two random variables $\xi_0$ and $\xi_1$ with two independent unbiased coin tosses : 
$\xi_0$ and $\xi_1$ are independent, and each one is equal to 1 with probability $1/2$. Then, we set
	$$ \xi_2\ \egdef\ \xi_0+\xi_1\mod 2. $$
Each random variable $\xi_i$ is called a $0$-block, and the triple $(\xi_0,\xi_1,\xi_2)$ is called a \emph{1-block}. We pick the second 1-block $(\xi_3,\xi_4,\xi_5)$ in the same way as the first one, but independently. The third 1-block $(\xi_6,\xi_7,\xi_8)$ is now set to be the pointwise sum of the first two 1-blocks:
\begin{gather*}
\xi_6\ \egdef\ \xi_0+\xi_3\mod 2,\\
\xi_7\ \egdef\ \xi_1+\xi_4\mod 2,\\
\xi_8\ \egdef\ \xi_2+\xi_5\mod 2.
\end{gather*}
Observe that this third 1-block follows the same distribution as the first two: $\xi_6$ and $\xi_7$ are two independent Bernoulli random variables with parameter $1/2$, and $\xi_8$ is the sum mod~2 of these variables. Note also that the third 1-block is independent of the first one, independent of the second one, but of course not independent of the first two together. The 9-tuple $(\xi_0,\xi_1,\ldots,\xi_8)$ is called a \emph{2-block}.

We can repeat this procedure inductively to construct $k$-blocks for each $k\ge0$: Suppose that for some $k$ we already have constructed the first $k$-block, which is the $3^k$-tuple $(\xi_0,\ldots,\xi_{3^k-1})$. Then, choose the second $k$-block $(\xi_{3^k},\ldots,\xi_{2\times 3^k-1})$ with the same probability distribution, but independently of the first one, and set the third $k$-block $(\xi_{2\times 3^k},\ldots,\xi_{3^{k+1}-1})$ to be the pointwise sum of the first two $k$-blocks:
\begin{equation}
\label{block-3-dot-rule}
\xi_{2\times 3^k+j}\ \egdef\ \xi_j+\xi_{3^k+j}\mod 2\ (0\le j\le 3^k-1).
\end{equation} 

\begin{figure}[h]
\centering
\includegraphics{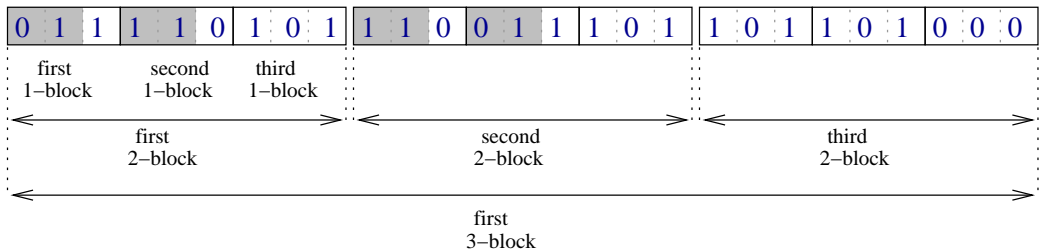}
\caption{Block construction of a stochastic process: The shaded coordinates are given by independent coin tosses . The non-shaded coordinates are computed from the shaded ones by 3-dot-type rules.}
\label{fig:blockConstruction}
\end{figure}

This inductive procedure gives the construction of a stochastic one-dimensional process $(\xi_i)_{i\ge0}$.
(This construction can easily be extended to a process indexed by $\ZZ$: Set the negative coordinates independently of the nonnegative ones by a similar symmetric construction.) Let us sketch the proof that our process is 2-fold mixing. For this, we use the two following facts, whose verification is left to the reader:
\begin{itemize}
\item Two different $k$-blocks are always independent.
\item Call a \emph{$k$-overlapping} the concatenation of two consecutive $(k-1)$-blocks lying in two different $k$-blocks. Any $k$-overlapping is independent of any concatenation of two consecutive $(k-1)$-blocks lying in other $k$-blocks.
\end{itemize}
Now, take two windows of fixed size $\ell$, and let $k$ be an integer such that $\ell\le 3^{k-1}$. Then, if the distance between the two windows is greater than $3^k$, either they lie in two different $k$-blocks, or at least one of them lie in a $k$-overlapping. In both cases the two windows are independent.

However, the stochastic process $\xi$ is clearly not 3-fold mixing, since for any $k\ge 0$, we have
$$ \xi_0 + \xi_{3^k} + \xi_{2\times 3^k}\ =\ 0 \mod 2. $$
This, of course, does not make $\xi$ a counterexample to Rohlin's question: The process we have just constructed is not a stationary one! Indeed, the pattern `111' for example can not be seen in the sequence $\xi_0\xi_1\xi_2$, but it can occur in the sequence $\xi_1\xi_2\xi_3$ with probability 1/8.

\subsection{How to make the construction stationary}
\label{stationarization}
The example described in the preceding section can be turned into a stationary process by applying some trick which is presented here. The process is still inductively constructed with $k$-blocks which follow the same distribution as before. The difference consists in the way $k$-blocks are extended to $(k+1)$-blocks.
Observe that a $k$-block lying in a given $(k+1)$-block can have three positions, which will be denoted by `0' (the first $k$-block in the $(k+1)$-block), `1' (the second one) and `2' (the third one). We are going to define the increasing family of $k$-blocks ($k\ge0$) containing the coordinate $\xi_0$ by using a sequence $S=(S_k)_{k\ge0}$ of independent, uniformly distributed random variables, taking their values in $\{0,1,2\}$. 

We start the construction by picking the first $0$-block $\xi_0$ in the usual way, with a coin toss. Now, we have to decide whether this $0$-block is in the first, second or third position in the $1$-block. This is done by using the first random variable $S_0$. Next, we complete the $1$-block by tossing a coin for the first missing variable, and setting the last one to be the sum mod~2 of the two others.
The extension from the $k$-block to the $(k+1)$-block containing $\xi_0$ goes on in a similar way: Once we have determined the $k$-block, we use the random variable $S_k$ to decide whether this $k$-block is in the first, second or third position in the $(k+1)$-block. Then the first missing $k$-block is chosen independently, and the last one is set to be the pointwise sum of the two other $k$-blocks.

\begin{figure}[h]
\centering
\input{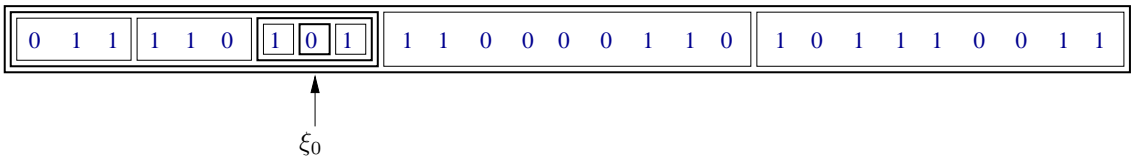}
\caption{Beginning of the construction with the skeleton sequence $S_0=1$, $S_1=2$, and $S_2=0$.}
\label{fig:skeleton}
\end{figure}

The embedding of $k$-blocks in $k+1$-blocks is called the \emph{skeleton} of the process, and the i.i.d. sequence $(S_k)_{k\ge0}$ coding this embedding is the \emph{skeleton sequence}. Since almost every realization of the skeleton sequence contains infinitely many 1's, the preceding procedure applied for all $k\ge0$ gives rise to $k$-blocks extending arbitrarily far away from 0 on both sides with probability one. This defines the whole process $\xi=(\xi_i)_{i\in\ZZ}$. 

Let us see how the skeleton sequence evolves when a coordinate shift is applied to the process $\xi$. It is not difficult to convince oneself that a shift of one coordinate to the left corresponds to the addition of `1' on the 3-adic number defined by the sequence $(S_k)$. (Write the sequence from right to left, and see it as a ``number'' written in base 3 with infinitely many digits, $S_0$ being the unitary digit; then add `1' to the sequence as you would do it for an ordinary number: Add `1' to $S_0$, and if $S_0$ reaches 3, then set $S_0=0$ and add `1' to $S_1$, and so on.) Observe that the distribution of the skeleton sequence is the same after this addition of `1', hence the distribution of the skeleton is invariant under the action of the coordinate shift. But once the skeleton is fixed, the distribution of the process is entirely determined by giving the distribution of $k$-blocks for every $k\ge0$, which is the distribution described in the preceding section. Therefore the whole distribution of the process is invariant under the coordinate shift, and the process we get now is stationary.

\begin{figure}[h]
\centering
\input{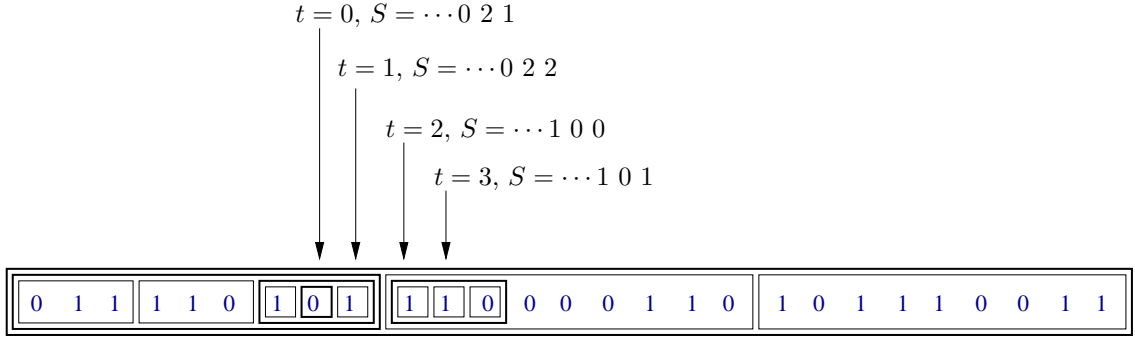}
\caption{Action of the coordinate shift on the skeleton sequence. The arrows denote the position of $\xi_0$ at successive times during the iteration of the shift.}
\label{fig:skeleton-shift}
\end{figure}

Unfortunately, making the process stationary has a cost: We have lost the 2-fold mixing property! Indeed, if for example we look at the rightmost coordinate $S_0$ of the skeleton sequence, we see that any realization of the process $\xi$ gives rise for $S_0$ to the periodic sequence 
$$ \cdots012012012\cdots $$
So, the process we get when we only observe $S_0$ is periodic. But if $\xi$ was 2-fold mixing, then every \emph{factor} of $\xi$ (that is to say, every stationary process which can be seen as a function of $\xi$, such as the process generated by $S_0$ for example\footnote{We leave as an exercise for the reader the verification of the fact that the skeleton sequence is indeed a function of $\xi$.}) would also be 2-fold mixing.

\subsection{A relative counterexample to Rohlin's question}

The stationary process generated by the whole skeleton sequence $S$ is well-known in ergodic theory, and is called the \emph{3-adic odometer}. (Be careful: $S$ does not take its values in a finite alphabet, it has infinitely many coordinates taking their values in $\{0,1,2\}$.) This process, which has appeared as a factor of $\xi$ in our new construction, is far from being mixing, since each of its coordinates is periodic. However, we can notice some interesting facts regarding the 2-fold and 3-fold mixing properties of $\xi$. Namely, once the skeleton is fixed (that is to say, conditionally to the $\sigma$-algebra generated by $S$), the mixing properties of $\xi$ are similar to those of the non-stationary process constructed in Section~\ref{non-stationary}. Thus, the process $\xi$ is 2-fold, but not 3-fold mixing \emph{relatively to the factor $\sigma$-algebra generated by $S$}. (More details on relative $k$-fold mixing can be found in \cite{ruedl13}.) 

It is a common idea in abstract ergodic theory to say that the study of stationary processes relatively to their factor $\sigma$-algebras gives rise to similar results as in the absolute study (one of the best examples of this fact is Thouvenot's relative version of Ornstein's isomorphism theorem \cite{thouv4}; another example is the proof of Proposition~\ref{zero-entropy-proposition} presented below). Therefore, the process that we have just constructed could make us think that a one-dimensional counterexample to Rohlin's question should exist. However, we are going to show in the next Section that, if such a process exists, it must be of a different nature than this one, or than Ledrappier's counterexample in $\ZZ^2$.

\section{Obstruction to the construction of a 3-dot-type one-dimensional counterexample}

\subsection{Multifold mixing and self-joinings}
\label{joining-section}
We need now to present a powerful tool which has been introduced in ergodic theory by Furstenberg~\cite{furst1}: The notion of self-joining of a stationary process. Let $\xi=(\xi_i)_{i\in\ZZ}$ be a stationary process taking its values in the alphabet $\AA$, and denote by $\mu$ its probability distribution on $\AA^{\ZZ}$. Take $\xi'$ another process defined on the same probability space, taking its values in the same alphabet $\AA$, and following the same distribution $\mu$. Then we can consider the joint process $(\xi,\xi')$, taking its values in the Cartesian square $\AA\times\AA$. If this joint process is still stationary, then we say that its distribution $\lambda$ on $\AA^{\ZZ}\times\AA^{\ZZ}\thickapprox(\AA\times\AA)^{\ZZ}$ is a \emph{2-fold self-joining of $\xi$}. In other words, a 2-fold self-joining of $\xi$ is a probability distribution on $\AA^{\ZZ}\times\AA^{\ZZ}$ whose marginals are both equal to $\mu$, and which is invariant under the coordinate shift.

Let us see some simple examples of such self-joinings. The first idea is to take the two processes $\xi$ and $\xi'$ independent of each other. Then we get $\mu\otimes\mu$ as our first example of a 2-fold self-joining. 
Another very simple example is obtained by taking $\xi'=\xi$, and we denote by $\Delta_0$ (``diagonal measure'') the 2-fold self-joining of $\xi$ concentrated on the diagonal of $\AA^{\ZZ}\times\AA^{\ZZ}$. This can be generalized by considering the case where $\xi'$ is equal to a shifted copy of $\xi$: We fix some $p\in \ZZ$, and we set $\xi'_i\egdef\xi_{i+p}$ for each $i\in\ZZ$. We denote by $\Delta_p$ the shifted diagonal measure obtained in this way. 

The set $J_2(\xi)$ of all 2-fold self-joinings of $\xi$ is endowed with the metrizable topology defined by the following distance:
$$
d(\lambda_1,\lambda_2)\ \egdef\ \sum_{n\ge0}\sum_{n'\ge0}\dfrac{1}{2^{n+n'}}
	\left| \lambda_1(\xi\in C_n,\xi'\in C_{n'}) - \lambda_2(\xi\in C_n,\xi'\in C_{n'}) \right|,
$$
where $(C_n)_{n\ge0}$ is the countable collection af all cylinder sets in $\AA ^{\ZZ}$. This topology (which is nothing else than the weak topology restricted to the set of 2-fold self-joinings of $\xi$) turns $J_2(\xi)$ into a compact metrizable topological space. The link with the 2-fold mixing property is now straightforward: The stationary process $\xi$ is 2-fold mixing if and only if the sequence $(\Delta_p)$ of shifted diagonal measures converges in $J_2(\xi)$ to the product measure $\mu\otimes\mu$ as $p\to +\infty$.

To translate the 3-fold mixing property into the language of joinings, we have to generalize the notion of self-joining to the case where 3 processes $\xi$, $\xi'$ and $\xi''$ with the same distribution $\mu$ are involved. This naturally leads to the definition of a \emph{3-fold self-joining of $\xi$}. (We can of course define an $r$-fold self-joining of $\xi$ for any $r\ge2$, but for our purpose the cases $r=2$ and $r=3$ will suffice.) The set $J_3(\xi)$ is also turned into a compact metrizable space when endowed with the restriction of the weak topology. Particularly simple and interesting elements of $J_3(\xi)$ are again the product measure $\mu\otimes\mu\otimes\mu$ and the shifted diagonal measures $\Delta_{p,q}$, $p,q\in\ZZ$, the latter denoting the distribution of the triple $(\xi,\xi',\xi'')$ when for all $i\in\ZZ$,
\begin{equation}
\label{Delta_pq}
\xi'_{i}=\xi_{i+p}\quad\mbox{and}\quad \xi''_{i}=\xi_{i+p+q}.
\end{equation}

The process $\xi$ is 3-fold mixing if and only if the following convergence holds in $J_3(\xi)$:
\begin{equation}
\Delta_{p,q}\ \tend{p,q}{+\infty}\ \mu\otimes\mu\otimes\mu.
\end{equation}

Now, let us assume that $\xi$ is a 2-fold mixing stationary process which is not 3-fold mixing. Then, since $\Delta_{p,q}$ does not converge to the product measure, we can find a subsequence $\Delta_{p_n,q_n}$ converging to some 3-fold self-joining $\lambda\neq\mu\otimes\mu\otimes\mu$. But the 2-fold mixing property of $\xi$ tells us that, under $\lambda$, the 3 processes $\xi$, $\xi'$ and $\xi''$ have to be pairwise independent. Hence, we get the following conclusion:

\begin{prop}
\label{pij}
If $\xi$ is a 2-fold mixing stationary process which is not 3-fold mixing, then $\xi$ has a 3-fold self-joining $\lambda\neq\mu\otimes\mu\otimes\mu$ with pairwise independent coordinates.
\end{prop}

\subsection{Restriction to zero-entropy processes}
\label{zero-entropy}
The natural question now is whether stationary processes satisfying the conclusion of Proposition~\ref{pij} can exist. But, without extra requirements, it is easy to find examples of such pairwise independent self-joinings which are not the product measure: Let $\xi$ consist of i.i.d. random variables $(\xi_p)_{p\in\ZZ}$, taking their values in $\{0,1\}$, each one with probability 1/2. Take an independent copy $(\xi'_p)_{p\in\ZZ}$ of this process and set $$ \xi''_p\ \egdef\ \xi_p+\xi'_p\quad \mod 2. $$
Then the three processes $\xi$, $\xi'$ and $\xi''$ have the same distribution, are pairwise independent,
but the 3-fold self-joining of $\xi$ we get in this way is not the product measure.

However, this process is not a counterexample to Rohlin's question. Indeed, all its coordinates being independent, $\xi$ is of course $k$-fold mixing for any $k\ge2$. It was Thouvenot who explained that this kind of Bernoulli shift situation could be avoided when one studies Rohlin's question, because it is always possible to restrict the analysis to the case of zero-entropy stationary processes. A stationary process $\xi=(\xi_i)_{i\in\ZZ}$ taking its values in a finite alphabet has zero entropy if and only if $\xi_0$ is measurable with respect to the \emph{past}, \textit{i.e.} with respect to the $\sigma$-algebra $\sigma(\xi_i,\ i<0)$. (Main definitions and results concerning entropy in ergodic theory can be found \textit{e.g.} in \cite{cornf1,glasn4,kalik5}.) 

\begin{prop}
\label{zero-entropy-proposition}
If there exists a stationary process $\xi$ which is 2-fold, but not 3-fold mixing, then one can find such a counterexample in the class of zero-entropy stationary processes.
\end{prop}

\begin{proof}
One of the main ingredient to prove this result is the so-called \emph{Pinsker $\sigma$-algebra} (or \emph{tail field}) of the process, which is the $\sigma$-algebra
$$ \Pi(\xi)\ \egdef\ \bigcap_{p\in\ZZ} \sigma(\xi_{-\infty}^p)\ =\ \bigcap_{p\in\ZZ} \sigma(\xi_p^{+\infty}). $$
The Pinsker $\sigma$-algebra of $\xi$ is invariant by the coordinate shift $T$: $(\xi_n)\mapsto(\tilde\xi_n)$, where $\tilde\xi_n\egdef \xi_{n+1}$. Therefore, if we choose any $\Pi(\xi)$-measurable random variable $\zeta_0$, taking its value in a finite alphabet $\BB$, and if we set for all $i\in\ZZ$ 
$$ \zeta_i\ \egdef\ \zeta_0\circ T^i, $$
then the whole stationary process $\zeta=(\zeta_i)_{i\in\ZZ}$ is $\Pi(\xi)$-measurable. What is remarkable is that such a process always has zero entropy, and that any stationary process $\zeta$ with zero entropy which is a factor of $\xi$ is automatically $\Pi(\xi)$-measurable (see \textit{e.g.} \cite{glasn4}, Theorem~18.6). Besides, Krieger's finite generator theorem ensures that it is always possible to find such a factor process $\zeta$, taking its values in an alphabet $\BB$ containing only two letters, and generating the whole Pinsker $\sigma$-algebra (nice proofs of Krieger's theorem can be found in \cite{glasn4}, \cite{kalik5}, or \cite{rudol7}). We henceforth fix such a process $\zeta$ with 
$$ \Pi(\xi)=\sigma(\zeta_i, i\in\ZZ). $$
As a factor of $\xi$, $\zeta$ automatically inherits the 2-fold mixing property. It remains to show that if $\zeta$ is 3-fold mixing, then so is $\xi$. 

For this, let us fix three cylinder events $E_1$, $E_2$ and $E_3$, which are measurable with respect to $\xi_{0}^{\ell}$ for some $\ell\ge0$. We have to compute the limit, as $p$ and $q$ go to $+\infty$, of the quantity
\begin{equation*}
\begin{split}
& \PP\left( \xi_0^{\ell}\in E_1,\xi_{p}^{p+\ell}\in E_2,\xi_{p+q}^{p+q+\ell}\in E_3 \right) \\
& = \esp{ \espc{\ind{\xi_0^{\ell}\in E_1}}{\sigma(\xi_p^{+\infty})} \ind{\xi_{p}^{p+\ell}\in E_2} 
    \ind{\xi_{p+q}^{p+q+\ell}\in E_3} } \\
& = \esp{ \left( \espc{\ind{\xi_0^{\ell}\in E_1}}{\sigma(\xi_p^{+\infty})} - \espc{\ind{\xi_0^{\ell}\in E_1}}{\Pi(\xi)} \right) \ind{\xi_{p}^{p+\ell}\in E_2} 
    \ind{\xi_{p+q}^{p+q+\ell}\in E_3} } \\
& \quad + \esp{ \espc{\ind{\xi_0^{\ell}\in E_1}}{\Pi(\xi)} \ind{\xi_{p}^{p+\ell}\in E_2} 
    \ind{\xi_{p+q}^{p+q+\ell}\in E_3} }.
\end{split}
\end{equation*}
The martingale convergence theorem gives 
$$ \espc{\ind{\xi_0^{\ell}\in E_1}}{\sigma(\xi_p^{+\infty})}\ \xrightarrow[p\to +\infty]{L^1}\ \espc{\ind{\xi_0^{\ell}\in E_1}}{\Pi(\xi)}, $$
hence the first term can be bounded by $\varepsilon$ if $p$ is large enough. We are left with
\begin{equation*}
\begin{split}
& \esp{ \espc{\ind{\xi_0^{\ell}\in E_1}}{\Pi(\xi)} \ind{\xi_{p}^{p+\ell}\in E_2} 
    \ind{\xi_{p+q}^{p+q+\ell}\in E_3} } \\
& = \esp{ \espc{\ind{\xi_{-p}^{-p+\ell}\in E_1}}{\Pi(\xi)} \espc{\ind{\xi_{0}^{\ell}\in E_2}}{\sigma(\xi_q^{+\infty})}
    \ind{\xi_{q}^{q+\ell}\in E_3} },
\end{split}
\end{equation*}
and again the martingale convergence theorem allows us to replace $\espc{\ind{\xi_{0}^{\ell}\in E_2}}{\sigma(\xi_q^{+\infty})}$ with $\espc{\ind{\xi_{0}^{\ell}\in E_2}}{\Pi(\xi)}$ if $q$ is large enough.
We thus have proven
\begin{multline}
\label{conditional_3-fold_mixing}
 \PP\left( \xi_0^{\ell}\in E_1,\xi_{p}^{p+\ell}\in E_2,\xi_{p+q}^{p+q+\ell}\in E_3 \right)  \\
 -  \esp{ \espc{\ind{\xi_{0}^{\ell}\in E_1}}{\Pi(\xi)}\ \espc{\ind{\xi_{p}^{p+\ell}\in E_2}}{\Pi(\xi)}
   \ \espc{\ind{\xi_{p+q}^{p+q+\ell}\in E_3}}{\Pi(\xi)} }\ \tend{p,q}{+\infty}\ 0.
\end{multline}
Note that in this equation, the expectation can also be written as
$$ \esp{ \espc{\ind{\xi_{0}^{\ell}\in E_1}}{\Pi(\xi)}\ \espc{\ind{\xi_{0}^{\ell}\in E_2}}{\Pi(\xi)}\circ T^{p}\ \espc{\ind{\xi_{0}^{\ell}\in E_3}}{\Pi(\xi)}\!\circ\!T^{p+q}}. $$
Now, if the process $\zeta$ generating $\Pi(\xi)$ is 3-fold mixing, we get that this expectation converges, as $p$ and $q$ go to $+\infty$, to the product
\begin{equation*}
\begin{split}
& \esp{ \espc{\ind{\xi_{0}^{\ell}\in E_1}}{\Pi(\xi)} }\ \esp {\espc{\ind{\xi_{0}^{\ell}\in E_2}}{\Pi(\xi)} }\ 
\esp {\espc{\ind{\xi_{0}^{\ell}\in E_3}}{\Pi(\xi)}} \\
& = \PP\left( \xi_0^{\ell}\in E_1\right)\PP\left( \xi_{0}^{\ell}\in E_2 \right) \PP\left( \xi_{0}^{\ell}\in E_3 \right),
\end{split}
\end{equation*}
which means that $\xi$ is also 3-fold mixing.
\end{proof}

Let us make some comment about this proof. The main argument is a nice illustration of the principle that the behaviour of a stationary process conditionally to one of its factors gives rise to similar results as in the non-conditioned case. The key notion here is the so-called \emph{$K$-property}: We say that the stationary process $\xi$ has the $K$-property if its Pinsker $\sigma$-algebra $\Pi(\xi)$ is trivial. It is easy to show that this $K$-property implies 3-fold mixing: Just write the preceding proof until \eqref{conditional_3-fold_mixing} in the case where $\Pi(\xi)$ is trivial, and you get the result.  (In fact, the same argument gives that the $K$-property implies $k$-fold mixing for any $k$). Now, when a factor $\zeta$ of $\xi$ is given, we can also define the $K$-property of $\xi$ \emph{relatively to $\zeta$} (see \textit{e.g.} \cite{ruedl8}), and check that $\xi$ always has the $K$-property relatively to its Pinsker $\sigma$-algebra: This comes from the fact that any factor of $\xi$ with zero entropy is $\Pi(\xi)$-measurable. But this in turn gives that $\xi$ is 3-fold mixing \emph{relatively to $\Pi(\xi)$}: This is more or less what \eqref{conditional_3-fold_mixing} says. The end of the proof consists in checking that, if $\zeta$ is 3-fold mixing, and if $\xi$ is 3-fold mixing relatively to $\zeta$, then $\xi$ is 3-fold mixing.

\medskip

Now that we have reduced Rohlin's problem to the case of zero-entropy processes, we can ask the question of pairwise independent self-joinings in this zero-entropy class.

\begin{question}
\label{pairwise-independent-question}
Does there exist a zero-entropy, 2-fold mixing stationary process $\xi$, and a 3-fold self-joining $\lambda$ of $\xi$
for which the coordinates are pairwise independent but which is different from the product measure?
\end{question}

Note also that if the assumption of 2-fold mixing is dropped, we can again find some counterexample: Take a periodic stationary process $\xi$ taking its values in $\{0,1,2\}$, where $\xi_0$ is uniformly distributed in the alphabet, and $\xi_{i+1}=\xi_i+1\mod3$ for all $i\in\ZZ$. Take an independent copy $\xi'$ of $\xi$, and set 
$$ \xi''_i\ \egdef\ 2\xi'_i-\xi_i\mod 3. $$ 
Then $\xi$, $\xi'$ and $\xi''$ share the same zero-entropy distribution, they are pairwise independent but the 3-fold self-joining they define is not the product measure.

\subsection{3-dot-type pairwise independent self-joinings}

Let us consider now the two-dimen\-sional example of Ledrappier from the point of view of self-joinings. The definition of the shifted diagonal 3-fold self-joining $\Delta_{p,q}$ is formally the same as in~\eqref{Delta_pq}, but in this case $p$ and $q$ are both elements of $\ZZ^2$. Particularly interesting is the sequence $\Delta_{p_n,q_n}$, where for all $n\ge 0$,
$$ p_n\ \egdef\ (0,2^n)\quad\mbox{and}\quad q_n\ \egdef\ (2^n,0). $$
From the relation \eqref{3-dot-n}, we see that the sequence $\Delta_{p_n,q_n}$ converges as $n\to+\infty$ to the 3-fold self-joining $\lambda$ of $\xi$ under which the three coordinates $\xi$, $\xi'$ and $\xi''$ are pairwise independent (because $\xi$ is 2-fold mixing), but which is not the product measure since, for all $i\in\ZZ^2$,
\begin{equation}
\label{3-dot-type}
\xi''_i\ =\ \xi_i+\xi'_i \mod 2.
\end{equation}

Observe also that, for the non-stationary process constructed in Section~\ref{non-stationary}, the sequence of 3-fold self-joinings $\Delta_{3^n,2\times 3^n}$ also converges to a joining with pairwise independent marginals which, by~\eqref{block-3-dot-rule}, satisfies a relation similar to \eqref{3-dot-type}. Since this kind of self-joining seems to appear naturally when one tries to construct examples which are 2-fold but not 3-fold mixing, we introduce the (slightly more general) following definition.

\begin{definition}
Let $\xi$ be a stationary process taking its values in a finite alphabet $\AA$.
We call \emph{3-dot-type self-joining} a 3-fold self-joining $\lambda$ of $\xi$ which has pairwise independent marginals, and for which there exists a map $f:\ \AA\times\AA\longrightarrow\AA$ such that
\begin{equation}
\label{general-3-dot-type}
\forall i\in\ZZ,\ \xi''_i\ =\ f(\xi_i,\xi'_i)\quad(\lambda\mbox{-a.s.}).
\end{equation}
\end{definition}

The result that we shall prove now shows that there is no hope to find a 2-fold but not 3-fold zero-entropy stationary process which admits a 3-dot-type self-joining: The only 3-dot-type self-joinings which can be seen in zero entropy are those arising from \emph{periodic} (therefore non 2-fold mixing) processes, like in the example presented at the end of Section~\ref{zero-entropy}.
\begin{theo}
\label{MainTheorem}
If the stationary process $\xi$ admits a 3-dot-type self-joining, then
\begin{itemize}
\item either $\xi$ is a periodic process,
\item or $\xi$ has entropy at least $\log 2$.
\end{itemize}
\end{theo}

The proof of the theorem is based on the following lemma.

\begin{lemma}
\label{MainLemma}
Let $X$, $Y$ and $Z$ be 3 random variables, sharing the same distribution on the finite alphabet $\AA$. Assume that these random variables are pairwise independent, and that there exists a map $f:\ \AA\times\AA\longrightarrow\AA$ such that
\begin{equation}
\label{map}
Z\ =\ f(X,Y)\quad(\mbox{a.s.}).
\end{equation}
Then their common distribution is the uniform distribution on a subset of $\AA$ 
\end{lemma}

\begin{proof}
Taking a subset of $\AA$ instead of $\AA$ if necessary, we can assume that each letter of $\AA$ is seen with positive probability. Let us fix some $y\in \AA$, and condition with respect to the event $(Y=y)$. Since $X$, $Y$ and $Z$ are pairwise independent, we get for any $x,z\in\AA$
$$ \PP(X=x|Y=y)\ =\ \PP(X=x)\quad\mbox{and}\quad\PP(Z=z|Y=y)\ =\ \PP(Z=z). $$
But knowing $(Y=y)$, we have the equivalence
$$ (Z=z)\ \Longleftrightarrow\ (f(X,y)=z). $$
Since $X$ and $Z$ can take the same number of values, we deduce that for any $z\in\AA$, there exists
a unique $x$ such that 
$$ f(x,y)\ =\ z, $$
and that furthemore, this $x$ satisfies
\begin{equation}
\label{equal-probability}
 \PP(X=x)\ =\ \PP(Z=z).
\end{equation}
Finally, observe that we can condition on $(Y=y)$ for any $y\in\AA$, and therefore that \eqref{equal-probability} holds for any $x$ and $z$ in $\AA$ for which we can find a $y\in\AA$ satisfying $f(x,y)=z$. But since $X$ and $Z$ are independent, this is true for \emph{any} $x,z\in\AA$.
\end{proof}

\begin{proof}[Proof of Theorem~\ref{MainTheorem}]
Fix some $m\ge1$, and apply Lemma~\ref{MainLemma} with $X\egdef\xi_0^{m-1}$, $Y\egdef{\xi'}_{\!0}^{\,m-1}$ and $Z\egdef{\xi''}_{\!0}^{\,m-1}$ under the 3-dot-type self-joining of $\xi$: We get that $\xi_0^{m-1}$ is uniformly distributed on the subset of sequences in $\AA^m$ which are seen with positive probability.
Denote by $p_m$ the number of such sequences. The same argument applied with $m+1$ in place of $m$ gives that each of the $p_{m+1}$ possible sequences of length $m+1$ is seen with probability $1/p_{m+1}$. Therefore, any possible sequence of length $m$ has exactly $a_m\egdef p_{m+1}/p_m$ different ways to extend to some possible sequence of length $m+1$, and each of the possible extensions has a conditional probability $1/a_m$. Now, observe that since $\xi$ is stationary, the number of possible extensions of a sequence of length $m+1$ can not be greater than the number of possible extensions of the  sequence of length $m$ obtained by removing the first letter, which gives the inequality $a_{m+1}\le a_m$.
Therefore, there are two cases:
\begin{itemize}
\item Either $a_m=1$ for $m$ large enough, and in this case the process $\xi$ is 
periodic;
\item Or $a_m$ is always greater than or equal to 2, and then the entropy of  $\xi$ is at least $\log 2$. 
\end{itemize}
\end{proof}

\subsection{Difference between dimension one and two}

It seems interesting to analyze the result proved in the preceding section, to see precisely which obstacle prevents us from constructing a 3-dot-type counterexample in dimension one, although it is allowed in dimension two. At which point does the one-dimensional argument developped above stop applying to the case of dimension two? 

First, note that the reduction of Rohlin's question to zero-entropy processes is still valid in the case of two-dimensional random fields: We can still consider in this case the Pinsker $\sigma$-algebra $\Pi(\xi)$, with respect to which every zero-entropy factor of $\xi$ is measurable, and prove that if $\xi$ is 2-fold but not 3-fold mixing, then this property comes from $\Pi(\xi)$. (For some presentations of the Pinsker $\sigma$-algebra in the multidimensional case, we refer the reader to \cite{conze1, glasn3, kamin1, kamin2}. For the generalization of Krieger's finite generator theorem to $\ZZ^d$-actions, see \textit{e.g.} \cite{rosen1, danil1}.) As far as Ledrappier's 3-dot example is concerned, there is no need to take a factor since this stationary random field already has zero entropy. 

Next, Lemma~\ref{MainLemma} gives that if $\xi$ is a finite-valued stationary random field indexed by $\ZZ^2$ admitting a 3-dot-type self joining, then the distribution of $\xi$ on any finite window is uniform on all the configurations which are seen with positive probability. And indeed, we can easily check that this is the case for Ledrappier's construction
(for example, any allowed configuration in a rectangular $\ell_1\times \ell_2$ window has probability $2^{-(\ell_1+\ell_2-1)}$). And this seems to be the point of the argument where there is a difference between $\ZZ$ and $\ZZ^2$: Although this uniform-probability property implies periodicity for one-dimensional zero-entropy stationary processes, this is no more true for two-dimensional stationary random fields, as we can see with Ledrappier's example.

We can also observe that in the one-dimensional case, if the stationarity property is dropped the argument fails at the same point: The relative counterexample constructed in Section~\ref{stationarization}, that is to say the process conditioned on the skeleton sequence, also has the uniform-probability property (which means in the non-stationary setting that for any $i\le j$, all the possible configurations for the sequence $\xi_i^j$ have the same probability). However, this process has zero entropy: There exists almost surely some $k$ such that $S_k=2$, therefore $\xi_0$ lies in the third $k$-block in its $(k+1)$-block. Then $\xi_0$ can be computed from the first two $k$-blocks, which are measurable with respect to the past of the process.

\section{Further questions}

Many important results have already been presented around Rohlin's multifold mixing problem. Most of them consider some special category of stationary processes (\textit{e.g.} finite-rank processes, or processes with singular spectrum), and prove that in this category a pairwise-independent 3-fold self-joining has to be the product measure. 

We hope that the work presented here can be the beginning of a slightly different approach: Consider a stationary process which admits some kind of pairwise-independent 3-fold self-joining which is not the product measure, and see which other properties on the process this assumption entails. In this direction, the most natural generalization of the present study should be the following situation: Take a stationary process $\xi$ which admits a pairwise-independent 3-fold self-joining $\lambda$, and assume that the three coordinates $\xi$, $\xi'$ and $\xi''$ of this self-joining satisfy
\begin{equation}
\label{GPIJ}
\xi''\ =\ \varphi(\xi,\xi')\quad\lambda\mbox{-a.s.}
\end{equation}
for some measurable function $\varphi:\,\AA^\ZZ\times\AA^\ZZ\to\AA^\ZZ$. Then what can be said on the process $\xi$? Can it have zero entropy and be 2-fold mixing? 

Since this general question seems to be quite difficult, some more restricted classes of pairwise independent self-joinings could be considered first, satisfying~\eqref{GPIJ} with some regularity assumption on $\varphi$. For example, what happens if under $\lambda$, the coordinate $\xi''_0$ is a function of finitely many coordinates of $\xi$ and $\xi'$?


\providecommand{\bysame}{\leavevmode\hbox to3em{\hrulefill}\thinspace}
\providecommand{\MR}{\relax\ifhmode\unskip\space\fi MR }
\providecommand{\MRhref}[2]{%
  \href{http://www.ams.org/mathscinet-getitem?mr=#1}{#2}
}
\providecommand{\href}[2]{#2}

\end{document}